\documentclass[12pt]{amsart}    
\usepackage{amscd}

\let\To=\longrightarrow
\let\Union=\bigcup

\title{Monomial ideals with linear quotients 
whose Taylor resolutions are  minimal}
\author{Munetaka Okudaira  and Yukihide Takayama}

\address{Munetaka Okudaira, Graduate School of Mathematical
Sciences, Ritsumeikan University, 
1-1-1 Nojihigashi, Kusatsu, Shiga 525-8577, Japan}
\email{rp001012@se.ritsumei.ac.jp}

\address{Yukihide Takayama, Department of Mathematical
Sciences, Ritsumeikan University, 
1-1-1 Nojihigashi, Kusatsu, Shiga 525-8577, Japan}
\email{takayama@se.ritsumei.ac.jp}

\theoremstyle{plain}
\newtheorem{thm}{Theorem}[section]
\newtheorem{prop}[thm]{Proposition}
\newtheorem{lem}{Lemma}
\newtheorem{cor}{Corollary}
\theoremstyle{definition}

\newtheorem{exmp}{Example}
\theoremstyle{remark}

\newcommand{\set}[1]{{\rm set}(#1)}

\newcommand{\lcm}{{\rm lcm}}
\begin{document}

\maketitle

\begin{abstract}
We study when Taylor resolutions of monomial ideals 
are minimal, particularly for ideals with linear quotients.
\end{abstract}

\section*{Introduction}

Let $S = K[X_1,\ldots, X_n]$ be a polynomial ring
over a field $K$ and consider a monomial ideal 
$I \subset S$. 
Let  $G(I) = \{u_1,\ldots, u_r\}$ be the minimal set of monomial
generators of $I$. Then  the Taylor resolution $(T_\bullet(I),
d_\bullet)$
of $I$ is defined as follows (cf. \cite{E} Excer.~17.11): 
$T_q(I) = \bigwedge^{q+1}L$ for $q=0,\ldots, {r-1}$
where $L$ is the $S$-free module with the basis $\{e_1,\ldots, e_{r}\}$
and $d_q : T_{q}(I) \To T_{q-1}(I)$, for $q=1,\ldots, r-1$, is defined by 
\begin{equation*}
 d_q(e_{i_0}\wedge\cdots\wedge e_{i_q})
 = \sum_{s=1}^{q}(-1)^s
     \displaystyle{\frac{\lcm(u_{i_0},\ldots, u_{i_q})}
                        {\lcm(u_{i_0},\ldots,\check{u}_{i_s},\ldots
			u_{i_q})}}
    e_{i_0}\wedge\cdots \overset{i_s}{\check{\hbox{}}}\cdots \wedge e_{i_q}
\end{equation*}
and the augumentation $\epsilon : T_{0}\To I$ is defined by
$\epsilon(e_i) = u_i$ for $i=1,\ldots, r$. 
In general, $T_\bullet(I)$ is far from minimal and the aim
of this paper is to determine some of the cases in which this resolution
is minimal.

A monomial ideal $I \subset S$ is said to be an {\em ideal with linear
quotients} if, for some specified order $u_1,\ldots, u_r$ of the minimal
set of generators, the colon ideals $(u_1, \dotsc, u_{j-1}) : u_j$ are 
generated by a subset of $\{ X_1, \dotsc, X_n \}$, for $j=1,\ldots, r$.
When we consider such an ideal $I = (u_1,\ldots, u_r)$, we will always
assume that $I$ has linear quotients with this order of the minimal
set of generators $u_1,\ldots, u_r$.
We also set $\set{u_j} = \{i_1,\ldots, i_s\}$ when 
$(u_1,\ldots, u_{j-1}): u_j = (X_{i_1},\ldots, X_{i_s})$. 
Stable ideals, squarefree stable ideals and  (poly)matroidal ideals are all
ideals with linear quotients and they have Eliahou-Kervaire type 
minimal resolutions \cite{HT}.

We will show that an ideal $I = (u_1,\ldots, u_r)$ with linear quotients 
has a minimal Taylor resolution if and only if $\vert\set{u_i}\vert = i-1$ for 
$i =1,\ldots, r$ (Theorem~\ref{thm:LinQuoSet}), where $\vert A \vert$ 
denotes the cardinarity of the set $A$. In the case of stable ideal,
this is precisely when 
$1\leq r \leq n$ and $u_i$ are in the form of 
$u_i = X_i(\prod_{k=1}^{i}X_k^{n_k})$, $i=1,\ldots, r$, for some integers 
$n_1,\ldots, n_r\geq 0$ (Theorem~\ref{thm:main}).
On the other hand, for a monomial ideal $I\subset S$ with linear
resolution, $I$ has the minimal  Taylor resolution precisely when
$I$ is in the form of $I = u\cdot (X_{i_1}, \ldots, X_{i_k})$,
where $u$ is a monomial and $\{i_1,\ldots, i_k\}\subseteq \{1,\ldots,
n\}$ (Theorem~\ref{thm:linearcase}). Such an ideal also has 
linear quotients. We also give several examples such as  matroidal
ideals and squarefree stable ideal having minimal Taylor resolutions.

We thank J\"{u}rgen Herzog for valuable discussions and comments.

\section{Ideal with linear quotients}

This section recalls some general facts on ideal with linear quotients
and give a condition for such ideals to have the minimal Taylor
resolutions.

\begin{lem}\label{set}
Let $I = (u_1, \dotsc, u_r)$ be a monomial ideal with linear quotients. Then, 
\[ \vert\set{u_i}\vert \leq i-1 \qquad {\rm for} \ i = 1, \dotsc, r \]
\end{lem}
\begin{proof}
Since $\set{u_i} = \{j\;\vert\; X_j\in\Union_{k=1}^{i-1}(u_k):u_i\}
= \Union_{k=1}^{i-1}\{j \;\vert\; X_j\in (u_k):u_i\}$
and each $(u_k):u_i$ is generated by a single variable, we obtain the 
desired result.
\end{proof}

\begin{lem}[cf. lemma 1.5 \cite{HT}]\label{Herzog-Takayama}
\label{lem:ekform}
Let $I$ be a monomial ideal with linear quotients. Then the Betti 
numbers $\beta_q(I)$ of $I$ are as follows:
\begin{equation*}
\beta_q(I) = \sum_{u \in G(I)} \binom{\vert\set{u}\vert}{q}
\quad \mbox{for all }q\geq 0.
\end{equation*}
\end{lem}
Recall that a monomial ideal $I\subset S$ is stable 
if, for an arbitrary  monomial $w \in I$, we have 
$X_{i}w/X_{m(w)} \in I$ for all $i< m(w)$
where $m(u) = \max\{j \;\vert\; X_j\mbox{ divides } u\}$. 
If $I = (u_1,\ldots, u_r)$ is stable, we have $\set{u} = \{1,\ldots,
m(u)-1\}$ if $\deg u_1 \leq \cdots \leq \deg u_r$
and $u_i > u_{i+1}$ by reverse lexicographical order if 
$\deg u_i = \deg u_{i+1}$. Then we can recover
the well-known Eliahou-Kervaire formula \cite{EK} from Lemma~\ref{lem:ekform}.
\begin{thm}
\label{thm:LinQuoSet}
Let $I = (u_1,\ldots, u_r)$ be a monomial ideal with linear quotients.
Then $I$ has the minimal Taylor resolution
if and only if  $\vert\set{u_i}\vert = i-1$ for $i = 1, \dotsc, r$. 
\end{thm}
\begin{proof}
We have 
\begin{equation}
\label{1}
\beta_q(I)\leq \sum_{i=1}^{r} \binom{i-1}{q}
\quad\mbox{for all }q\geq 0
\end{equation}
by Lemma~\ref{set} and ~\ref{lem:ekform}. On the other hand, 
the Taylor resolution $T_\bullet(I)$ is minimal if and only if 
$\beta_q(I) = \binom{r}{q+1} = \sum_{i=1}^{r} \binom{i-1}{q}$. Thus 
the inequality in (\ref{1}) must be equality, which implies 
$\vert\set{u_i}\vert =i-1$ for all $i$ by Lemma~\ref{set}.
\end{proof}

\section{Stable ideals having the minimal Taylor resolutions}

The goal of this section is to determine precisely 
the stable ideals that have the minimal Taylor resolutions.

We fist prepare a formal characterization of such ideals.
\begin{prop}
\label{prop:minTaylorStable}
Let $I$ be a stable ideal of $R$. Then the following conditions are equivalent :
\begin{enumerate}
\item $I$ has the mimimal Taylor resolution;
\item $\max\{ m(u)\;\vert\;u \in G(I) \} = \vert G(I)\vert$ ;
\item $m_{i}(I) =
\begin{cases}
 1 & \text{for $1 \leq i \leq \vert G(I)\vert$} \\
 0 & \text{for $\vert G(I)\vert < i \leq n$},
\end{cases}$
\end{enumerate}
where we define $m_i(I) = \vert\{u\in G(I) \;\vert\;  m(u) = i  \}\vert$,
\end{prop}
\begin{proof} We first show that $m_i(I)\geq 1$ for $i=1,\ldots, b_0$,
where $b_0 := \max\{ m(u)\;\vert\;u \in G(I) \}$.
We have $m_{b_0}(I)\geq 1$ and 
there exists $v\in G(I)$ such that 
$v = wX_{b_0}^{\alpha}$ with $m(w) < b_0$  and $\alpha>0$.
Suppose that there exists $j$ such that $1 \leq j < b_0$ and
$m_j(I) = 0$.
Since $I$ is stable, we have $v' := w X_j^{\alpha} \in I$ and 
there exists $u\in G(I)$ that divides $v'$. 
Since $m_j(I)=0$, $u$ must divide $w$, which implies
$u\in G(I)$ divides $v\in G(I)$, a contradiction.
Thus we have $b_0 \leq \vert G(I)\vert$.

Now we show (i) $\Rightarrow$ (ii) : Assume that 
$I$ has the minimal Taylor resolution, i.e.,
$\beta_{q}(I) = \binom{\vert G(I)\vert}{q+1} = \sum_{i=1}^{\vert G(I)\vert}\binom{i-1}{q}$
for $q\geq 0$.
By Lemma~\ref{lem:ekform} and Theorem~\ref{thm:LinQuoSet}, we have $\beta_{q}(I) = 
\sum_{i=1}^{b_0}\binom{i-1}{q}$, so that we must have $b_0 = \vert G(I)\vert$.
(ii) $\Rightarrow$ (iii) is clear from the inequality $b_0 \leq \vert G(I)\vert$.
Finally we show (iii) $\Rightarrow$ (i) : 
By (iii), we have  $G(I) = \{u_1, \dotsc, u_r\}$ 
with $m(u_i) =i$ for $i=1,\ldots, r$. We claim that 
with this order of the generators $I$ is an ideal with linear
quotients with $\set{u_i} = \{1,\ldots, m(u_i)-1\}$.
To prove the claim we have only to show that
$\deg{u_i}\leq\deg{u_{i+1}}$
for all $i$ (see the comment after Lemma~\ref{lem:ekform}).
We can set $u_i = vX_i^p$ and $u_{i+1}=v'X_{i+1}^q$ for some monomials
$v$ and $v'$ with $m(v)< i$ and $m(v') < i+1$ and integers $p,q>0$
Then, since $I$ is stable, we have $w:= v'X_i^q \in I$
so that there exists $u_j\in G(I)$ that divides $w$. In particular,
$m(u_j)\leq i$ and $\deg{u_j}\leq \deg{w}=\deg{u_{i+1}}$. If $m(u_j)<i$, then $u_j$ 
divides $v'$, which implies that $u_j$ divides $u_{i+1}$, a
 contradiction.
Thus $m(u_j)=i$ so that by (iii) we must have $u_j=u_i$.
Then $\deg{u_j}\leq \deg{u_{i+1}}$ as required.
Now we have $\vert\set{u_i}\vert =i-1$, so that $I$ has the mimimal
Taylor resolution by Theorem~\ref{thm:LinQuoSet}.
\end{proof}
Using above proposition we can determine the stable ideals
with the  minimal Taylor resolutions.
\begin{thm}
\label{thm:main}
Let $I\subset S$ be a stable ideal. Then $I$ has the mimimal
Taylor resolution if and only if it is in the following form:
$I = (u_1,\ldots, u_r)$ for some $r\leq n$, where 
\begin{equation*}
u_i =  X_i\cdot{\prod_{k=1}^i X_{k}^{a_k}}
\quad (i=1,\ldots, r)
\end{equation*}
for some integers $a_{1}, \dotsc, a_{r} \geq 0$.
\end{thm}
\begin{proof}We can easily check that an ideal in the above form
is stable and it has the minimal Taylor resolution by
Prop.~\ref{prop:minTaylorStable}. 
Now we show the converse.
By Prop.~\ref{prop:minTaylorStable} we can assume that 
$G(I) = \{u_1,\ldots, u_r\}$ with $m(u_i)=i$ for $i=1,\ldots, r$
and $r\leq n$. 
Thus we can write $u_i$ as follows
\begin{equation*}
u_i =  X_i\cdot{\prod_{k=1}^i X_{k}^{a_{i,k}}} \quad (i=1,\ldots, r)
\end{equation*}
for some integers $a_{i,k}\geq 0$, $1\leq i\leq r$ and $1\leq k\leq i$.
We will show that each $a_{i,k}$ is constant with regard to $i$.

Since $I$ is stable, we have 
$w := u_rX_{r-1}/X_{r} =  X_{r-1}\cdot{\prod_{k=1}^r X_{k}^{a_{r,k}}} \in
 I$,
so that there exists $u\in G(I)$ that divides $w$. We claim that
$u = u_{r-1}$. In fact, we have $m(u)\leq r-1$ since $u\ne n_r$.
If $m(u) < r-1$, then $u$ divides $w/X_{r-1} =\prod_{k=1}^r X_{k}^{a_{r,k}}$,
which implies that $u$ divides $u_r$, contradicting the assumption 
that both $u$ and $u_r$ are from $G(I)$. Thus we must have $m(u) = r-1$,
which implies $u=u_{r-1}$. 
Then we know that $a_{r-1,r-1} \leq a_{r,r-1}$.
If $a_{r-1,r-1}< a_{r,r-1}$, then $u_{r-1}$ divides
$u_r$,  a contradiction. Thus we have  $a_{r-1,r-1}=a_{r,r-1}$.
Repeating the same argument replacing $u_r$ by $u_j$ 
$(j=r-1, r-2, \ldots, 2)$, we obtain
$a_{j,j}     =    a_{j+1,j}$ for $j=1,\ldots, r-1$.

Now assume by induction that we have already obtained that 
$u_j$ divides $w_j := u_r(X_j/X_r)
= X_j(\prod_{i=1}^rX_i^{a_{r,i}})$ 
for all $j$ such that $k\leq j \leq r-1$ and consider 
$w_{k-1}$. By the similar discussion to the above, there 
exists $u\in G(I)$ such that $u$ divdes $w_{k-1}$ and 
$m(u)\geq k-1$. Assume that $j_0:= m(u) > k-1$. We have 
$u= u_{j_0}$ and by assumption
$u$ divides both $X_{j_0}(\prod_{i=1}^rX_i^{a_{r,i}})$ and 
$X_{k-1}(\prod_{i=1}^rX_i^{a_{r,i}})$, which 
implies that $u$ divides $u_r$, a contradiction. Thus we have 
$u = u_{k-1}$. Consequently, $u_j$ divides 
$w_j$ for all $j=1,\ldots, r-1$ and by the similar discussion to
the above we have $a_{j,j}= a_{r, j}$ for all $j=1,\ldots, r-1$.

Carrying out the same discussion replacing $u_r$ by $u_\ell$
for $\ell = r-1, r-2,...$, we obtain 
$a_{j,j}= a_{i, j}$ for all $i=1,\ldots, r$ and $j=1,\ldots, r-1$,
as required.
\end{proof}

\begin{cor}\label{example1}
Assume that $I\subset S$ is a stable ideal generated by monomials 
of the same degree $d>1$. Then $I$ has the minimal Taylor resolution
if and only if $I$ is in the form of 
$X_1^{d-1}(X_1,X_2,\ldots, X_r)$ for some $r\leq n$.
\end{cor}

\section{Linear minimal Taylor resolutions}
In this section, we consider non-stable cases of ideals with linear quotients.

\begin{thm}
\label{thm:linearcase}
Let $I$ be a monomial ideal with a linear resoluition. Then the following conditions are equivalent:
\begin{enumerate}
\item $I$ has the mimimal Taylor resolution:
\item  $I = u \cdot (X_{i_1}, \dotsc, X_{i_k})$  for some 
$1\leq i_1 < \cdots < i_k \leq n$ and a monomial $u$.
\end{enumerate}
In this case, $I$ is an ideal with linear quotients.
\end{thm}
\begin{proof} It is clear that ideals in the form of (ii)
have linear quotients. 
We only have to show (i) to (ii), and the converse is clear.
We prove by induction on $\vert G(I)\vert$. 
Let $G(I) = \{u_1,\ldots, u_r\}$ $(r\geq 2)$
and let $T_\bullet(I)$ be the linear minimal Taylor resolution.
Since we have 
\begin{equation*}
\frac{\lcm(u_1, u_{i_1},\ldots, u_{i_q})}
     {\lcm(u_1, u_{i_1},\ldots,\check{u}_{i_s},\ldots u_{i_q})}
= \frac{\lcm(u_{i_1},\ldots, u_{i_q})}
     {\lcm(u_{i_1},\ldots,\check{u}_{i_s},\ldots u_{i_q})}
\end{equation*}
for all $1 < i_1 < \ldots, i_q\leq r$, by truncating 
all the bases in the form of $e_1\wedge\cdots$
from $T_\bullet(I)$,
we obtain the linear minimal Taylor resolution of $J=(u_2,\ldots, u_r)$.
Thus by the induction hypothesis we have 
$u_k = uX_{i_j}$, $k=2,\ldots, r$, for some monomial $u$ and 
$1\leq i_1<\cdots <i_j\leq n$. Now we show that $u_1$ is in
the form of $uX_{i_0}$ for some $i_0\notin \{i_1,\ldots, i_j\}$.
Since $T_\bullet(I)$ is linear, both $\lcm(u_1, u_i)/u_1$  
and $\lcm(u_1, u_i)/u_i$ must be linear for $i=2,\ldots, r$. 
Then we easily know that $u$ must divide $u_1$ and 
conclusion follows.
\end{proof}

Now we show some examples produced by 
Theorem~\ref{thm:linearcase}.

\begin{exmp}
Let $I=X_1^{p}X_{k+1}^q(X_1,X_2,\ldots, X_k)$ with $1<k\leq n$
and $p,q\geq 1$.
Then $I$ is a non-stable ideal with linear quotients whose Taylor
resolution is minimal.
\end{exmp}

\begin{exmp} A Stanley-Reisner ideal $I \subset S$ generated 
by squarefree monomials with the same degree is called 
matroidal if it satisfies the following exchange property:
For all $u, v\in G(I)$ and all $i$ with $\nu_i(u) > \nu_i(v)$,
there exists an integer $j$ with $\nu_j(u) < \nu_j(v)$ such that 
$X_j(u/X_i)\in G(I)$, where we define
$\nu_i(u) = a_i$ for $u = X_1^{a_1}\cdots X_n^{a_n}$.
A matroidal ideal $I$ has a linear resolution (cf. \cite{HT}).
If it has the minimal Taylor resolution, we know that $I$ 
is in the form of 
\[ I = X_{i_1} \dotsb X_{i_{p}} (X_{j_{1}}, \dotsc, X_{j_{q}})  \]
for $\{i_{1}, \dotsc, i_{p} \} \cap \{j_{1}, \dotsc, j_{q} \} = \emptyset$ 
with $p+q\leq n$.
\end{exmp}
\begin{exmp}
A  squarefree stable ideal is a Stanley-Reisner ideal $I\subset R$ 
satisfying the condition that, for all $i< m(u)$ such that $X_i$ does not divide
$u$, one has  $X_i(u/X_{m(u)})\in I$. 
Let $I$ be a squarefree stable ideal generated by monomials 
with the same degree.
If $I$ has the minimal Taylor resolution, then 
$I$ is in the form of $I = u(X_{p+1},\ldots, X_q)$ with
$u = X_1\cdots X_p$ for some $1\leq p\leq q$.
\end{exmp}

\end{document}